\documentclass[english,12pt]{amsart}
\usepackage[english]{babel}
\usepackage{amsmath,amssymb,enumerate,amsthm}
\usepackage[latin1]{inputenc}
\usepackage[T1]{fontenc}
\usepackage{cite}
\usepackage{latexsym}
\usepackage{setspace}

\newtheorem{lemma}{Lemma}[section]
\newtheorem{theorem}[lemma]{Theorem}
\newtheorem{prop}[lemma]{Proposition}

\newtheorem{rk}[lemma]{Remark}

\newtheorem*{p-adic}{$p$-adic Littlewood Conjecture}
\newtheorem*{Mixed}{Mixed Littlewood Conjecture}
\newtheorem{hypothesis}{Hypothesis}
\theoremstyle{remark}

\theoremstyle{plain}

\newcommand{\s}{\mathcal{S}}

\def\={\;=\;}
\def\>{\;>\;}
\def\<{\;<\;}
\def\:{\,:\,}
\def\.={\;\dot{=}\;}

\newcommand{\R}{\mathbb{R}}

\newcommand{\N}{\mathbb{N}}
\newcommand{\D}{\mathcal{D}}
\newcommand{\cO}{\mathcal{O}}

\newcommand{\nn}{t}

\begin{document}

\title[Mixed Littlewood Conjecture]{On the Mixed Littlewood Conjecture and continued fractions in quadratic fields}
\author{Paloma BENGOECHEA, Evgeniy ZORIN}
\address{Department of Mathematics, University of York, York, YO10 5DD, United Kingdom}
\email{paloma.bengoechea@york.ac.uk, evgeniy.zorin@york.ac.uk}

\begin{abstract}
We show how a recent result by Aka and Shapira on the evolution of continued fractions in a fixed quadratic field implies the classic result of de Mathan and Teulli\'e on the Mixed Littlewood Conjecture.
\end{abstract}

\maketitle

\section{Introduction}

A famous open problem in the field of Diophantine approximation is the Littlewood Conjecture which claims that, for every pair $(\alpha,\beta)$ of real numbers, we have
\begin{equation}\label{Lc}
\inf_{q\geq 1}\, q\cdot\|q\alpha\|\cdot\|q\beta\|\=0,
\end{equation}
where $\|\cdot\|$ denotes the distance to the nearest integer. 
The first significant contribution to this question is due to Cassels and Swinnerton-Dyer \cite{CS} who proved that \eqref{Lc} is satisfied when $\alpha$ and $\beta$ belong to the same cubic field. This result was sharpened by Peck \cite{P} who showed that if $1,\alpha,\beta$ form a basis of a cubic field, then 
$$
\liminf\limits_{q\rightarrow\infty}\, q\cdot \log q \cdot \|q\alpha\| \cdot \|q\beta\|\<\infty.
$$
Further examples of pairs $(\alpha,\beta)$ satisfying \eqref{Lc} are given in \cite{AB} and \cite{dM}.
Despite some recent progress on the Hausdorff dimension of the set of counterexamples  \cite{PV},\cite{EKL}, the Littlewood Conjecture remains an open problem. 



In 2004, de Mathan and Teulli\'e~\cite{dMT} proposed a variant of the Littlewood Conjecture, called Mixed Littlewood Conjecture,
in which the quantity $\|q\beta\|$ is replaced by a pseudo-absolute value $|q|_\mathcal{D}$. A pseudo-absolute sequence $\D$ is an increasing sequence of positive integers $\mathcal{D}=(u_n)_{n\in\mathbb{N}}$ with $u_1=1$ and $u_n | u_{n+1}$ for all $n$. The pseudo-absolute value $|q|_\mathcal{D}$ is then defined by 
$$
|q|_\mathcal{D}=\inf\left\{1/u_n\, :\, q\in u_n\mathbb{Z}\right\}.
$$
When $\mathcal{D}$ is the sequence $(p^n)_{n\in\N}$, where $p$ is a prime number, then $|\cdot|_\mathcal{D}$ is the usual $p$-adic value $|\cdot|_p$, normalised such that $|p|_p=p^{-1}$, and de Mathan and Teulli\'e's conjecture is then known as $p$-adic Littlewood Conjecture.

\begin{Mixed} For every real number $\alpha$ and every pseudo-absolute sequence $\D$, we have
\begin{equation}\label{p-adicc}
\inf_{q\geq 1}\, q\cdot\|q\alpha\|\cdot|q|_{\mathcal{D}}\=0.
\end{equation}
\end{Mixed}

This conjecture obviously holds when $\alpha$ is rational or has unbounded partial quotients. Thus one only has to consider the case when $\alpha$ is a badly approximable number (that is irrational $\alpha$ with bounded partial quotients). From a metric point of view, the set of badly approximable numbers is moderately small: it has Lebesgue measure 0 but Hausdorff dimension 1.

In 2007, Einsiedler and Kleinbock \cite{EK} established that, for every given prime number $p$, the
set of counterexamples $\alpha$ to the $p$-adic Littlewood Conjecture has Hausdorff dimension 0, so this set is much smaller than the set of all badly approximable numbers.
In an opposite direction, Badziahin and Velani \cite{BV} proved that, for every given pseudo-absolute sequence $\mathcal{D}$, the set of real numbers $\alpha$ satisfying
\begin{equation}\label{BVteo}
\inf_{q\leq 3}\, q \cdot \log q \cdot \log\log q\cdot \|q\alpha\| \cdot |q|_\D\>0
\end{equation}
has full Hausdorff dimension.

However, the first contribution to de Mathan and Teulli\'e's conjecture goes back to themselves in~\cite{dMT}.
A natural family of badly approximable numbers are the quadratic irrational numbers, since their continued fraction expansion is periodic. Using machinery from $p$-adic analysis, de Mathan and Teulli\'e proved
~\eqref{p-adicc} for any such $\alpha$ when $\D$ is bounded (i.e. the ratio of two consecutive terms of $\D$ is bounded). Moreover, for any quadratic irrational number $\alpha$ and any bounded pseudo-absolute sequence $\D$, let $\s$ be the set of positive integers defined by
\begin{equation} \label{def_s}
\s:=\left\{q\in\mathbb{N}\:  \|q\alpha\| \ll \dfrac{1}{q}\right\},
\end{equation}
they proved a stronger statement 
\begin{equation}\label{p-adic-log}
\liminf\limits_{q\in\s}\, q \cdot \log q \cdot \|q\alpha\| \cdot |q|_\D\<\infty.
\end{equation}
Later, in the same direction as \eqref{BVteo}, de Mathan \cite{dM} showed that this result can not be significantly improved: for any quadratic irrational number $\alpha$ there exists a positive constant $\lambda$ (depending only on $\alpha$) such that
\begin{equation}\label{p-adic-log2}
\liminf\limits_{q\in\s}\, q \cdot (\log q)^\lambda \cdot \|q\alpha\| \cdot |q|_\D\>0,
\end{equation}
where the set $\s$ is defined by~\eqref{def_s}. The common belief is that \eqref{p-adic-log2} holds for any $\lambda\geq 1$, also when the $\liminf$ is taken over the set $q\in\N$ of all positive integers, and not only on the set $q\in\s$. 

In this paper we give a new proof of \eqref{p-adic-log} by using a recent result of Aka and Shapira \cite{AS} on the evolution of the continued fraction expansions of $\nn\alpha$, where $\alpha$ is a fixed quadratic irrationality and $\nn$ belongs to the group of $M$-units $\cO_M^{\times}$ for a fixed set of primes $M$. Their result allows us to estimate accurately the geometric mean of the partial quotients in the period of $\nn\alpha$ in terms of the length of the period and a constant depending only on $\alpha$ and the set of prime divisors of $t$. Then we find, when $\alpha$ is a quadratic irrationality, an explicit subset $\s'$ of de Mathan's set $\s$ for which there exists a positive constant $c$ such that, for all $q\in\s'$, we have
\begin{equation}\label{p-adic-BZ}
 q \cdot \log q \cdot \|q\alpha\| \cdot |q|_\D\,\leq\, c.
\end{equation} 

There is a substantial difference between the Mixed Littlewood Conjecture and its restrictive bounded form: for every given bounded pseudo-absolute sequence, the set of all primes dividing one of the terms of the sequence is finite, whereas it can be infinite if the sequence is
 not bounded. This discrepancy between the general conjecture and its bounded form was already palpable in \cite{dMT}, as it was the obstruction to extend the proof of \eqref{p-adic-log} to the general Mixed Littlewood Conjecture. The same happens in this article again. Although Lemma \ref{lema} gives a very simple proof of the Mixed Litlewood Conjecture for irrational quadratics, for the stronger result \eqref{p-adic-BZ} we need the set of all primes dividing one of the terms of the given pseudo-absolute sequence to be finite. This condition is due to Aka-Shapira's result \cite{AS}, which we use more than once. 

The link between Aka-Shapira's result and the bounded Mixed Littlewood Conjecture is particularly interesting because their results and methods seem to be extendable to continued fractions which are not exactly periodic. This approach, if developed, could provide a uniform way to establish the bounded Mixed Littlewood Conjecture or the $p$-adic Littlewood Conjecture in a vast variety of cases. 

We assume that the reader is familiar with the classical results from the theory of continued fractions\footnote{Will the reader need a reference on this subject, we recommend him~\cite{C1957}.}.

\section{Main result and proofs}

Throughout the paper $\alpha$ denotes a quadratic irrational number (so the continued fraction expansion of $\alpha$ is periodic) and $\D=(u_n)_{n\in\N}$ denotes a sequence of positive integers such that $u_1=1$ and $u_{n} | u_{n+1}$ for all $n\in\N$.
For $n\geq 0$, we denote by $l^{(n)}$ the length of the period of the continued fraction expansion of $u_n\alpha$, by $a^{(n)}_k$ the $k$-th partial quotient of this continued fraction, and by $r^{(n)}_k$ the denominator of the $k$-th convergent of $u_n\alpha$. We denote by $c_\alpha\in(0,1)$ a Markov constant of $\alpha$, so
$$
\inf\limits_{q\geq 1}\, q\cdot\|q\alpha\|\geq c_\alpha.
$$
We use the following additional hypothesis on the sequence $\D=(u_n)_{n\in\N}$.
\begin{hypothesis} \label{H1}
Assume there exists a finite set of primes $M$ such that for all $n\in\mathbb{N}$, we have $u_n=\prod_{p\in M}p^{m_p}$, with $m_p\geq 0$ (in other terms, all the integers $u_n$ belong to the group of $M$-units $\cO_M^{\times}$). 
\end{hypothesis}

We define the set of positive integers
\begin{equation} \label{def_s_prime} 
\s'=\left\{q\in\mathbb{N}\: \exists n\geq 0\mid q=u_n r^{(n)}_{l^{(n)}-1}\right\}.
\end{equation}
For this concrete set $\s'$ we have the following result.
\begin{theorem}\label{BZ}
Assume 
Hypothesis~\ref{H1}.
Then, for every quadratic irrational number $\alpha$ there exists a constant $c>0$ such that 
\begin{equation}\label{ecT}
q \cdot \log q \cdot \|q\alpha\| \cdot |q|_{\D}\leq c
\end{equation}
for all $q\in\s'$. In particular,
\begin{equation}\label{ecT2}
\liminf\limits_{q\in\s'}\, q \cdot \log q \cdot \|q\alpha\| \cdot |q|_{\D}<\infty.
\end{equation}
\end{theorem}
\begin{rk}
The constant $c$ in the statement of Theorem~\ref{BZ} depends on the sequence $\D$. In fact, we will see in the proof that this constant depends only on the finite set of primes $M$ from Hypothesis~\ref{H1} rather than on the sequence $\D$ itself.
\end{rk}

In order to prove Theorem~\ref{BZ}, we need some information on the growth of the partial quotients in the period of $u_n\alpha$. All the results we need are stated in Lemma~\ref{lema} and Proposition~\ref{E}.

In Lemma~\ref{lema} we study the partial quotients of $t\alpha$, for $t\in\N$. We provide a uniform upper bound for all of them and a lower bound for the last partial quotient in the period of the continued fraction expansion of  $\nn\alpha$. This lower bound agrees with the upper bound up to a multiplicative constant independent of $t$, so the uniform upper bound in Lemma~\ref{lema} is the best possible up to a multiplicative constant which only depends on $\alpha$.

However, Proposition~\ref{E} ensures that, under Hypothesis 1, most of the partial quotients of $u_n\alpha$ are actually much smaller than the uniform upper bound from Lemma~\ref{lema} which depends on $\alpha$ and $n$. Indeed, we show in Proposition~\ref{E} that the geometric mean of the partial quotients in the period of $u_n\alpha$ is bounded above by an absolute constant independent of $n$. 

This proposition is analoguous to a result of Khintchine, which states that, for any subset $X$ of $\R$ with full measure, the geometric mean of the first $N$ partial quotients of any $x\in X$ converges to a constant $\gamma$ independent of $x$ as $N\rightarrow\infty$. The constant $\gamma$ was calculated later by L\'evy~\cite{L1937}.

Proposition~\ref{E} follows from a recent result of Aka and Shapira~\cite{AS} which, informally speaking, shows that under Hypothesis~\ref{H1} the continued fractions of the numbers $u_n\alpha$, where $\alpha$ is fixed, behave like a generic continued fraction as $n$ tends to infinity. This means that the results about statistics of the continued fractions of elements in a set of full measure (like Khintchine's) should have their analogues when we consider the continued fractions of the numbers $u_n\alpha$, with fixed $\alpha$ and $n\in\mathbb{N}$ sufficiently large. Aka and Shapira give the rate of convergence of the evolution of the continued fractions of $u_n\alpha$ toward the ``average'' behaviour as $n\rightarrow\infty$, as well as an error term which prevents us from transferring Khintchine-L\'evy's result in its exact form (the error term leads in Proposition~\ref{E} to a constant $\kappa$, which does not appear in Khintchine-L\'evy's result). 

\begin{lemma}\label{lema} Let $\alpha$ be a quadratic irrational number and let $\nn\geq 2$ be an integer. Then
\begin{enumerate}[(i)]
\item the partial quotients of $\nn \alpha$ are bounded above by $\dfrac{\nn}{c_\alpha}$, 

\item the last partial quotient of the period of $\nn\alpha$ is bounded below by $\lfloor\nn\alpha\rfloor$.
\end{enumerate}
\end{lemma}

\textbf{Proof.} From the theory of continued fractions, we have
\begin{equation}\label{ecL1}
r^{(1)}_k \|r^{(1)}_k \nn\alpha\|\<\dfrac{1}{a^{(1)}_{k+1}}\qquad(k\geq 0).
\end{equation}
Multiplying both sides of the inequality by $\nn$ we obtain
\begin{equation}\label{ecL2}
c_\alpha<\nn\, r^{(1)}_k \|r^{(1)}_k \nn\alpha\|\<\dfrac{\nn}{a^{(1)}_{k+1}}\qquad(k\geq 0),
\end{equation}
which proves (i).

In order to show (ii), we can assume that the continued fraction expansion of $\alpha$ is purely periodic, i.e. $\alpha> 1$ and its Galois conjugate $\overline{\alpha}$ belongs to the interval $(-1,0)$. Therefore the continued fraction expansion of the number $\beta=\nn\alpha-\lceil \nn\overline{\alpha}\rceil$ that we denote by $[\overline{b_0,b_1\dots,b_s}]$ is also purely periodic, and we have
\begin{equation} \label{cf_alpha}
\nn\alpha=\beta + \lceil \nn\overline{\alpha}\rceil = [b_0 + \lceil \nn\overline{\alpha}\rceil, \overline{b_1, \ldots, b_s, b_0}]
\end{equation}
with $b_0=\lfloor\beta\rfloor \geq \lfloor\nn\alpha\rfloor$.
\begin{flushright}
$\square$
\end{flushright}

\begin{rk} The Mixed Littlewood Conjecture for quadratic irrationals follows directly from Lemma \ref{lema}. Indeed, if a quadratic irrational number $\alpha$ is a counterexample to the conjecture, then all the partial quotients of $\left\{u_n\alpha\right\}_{n\in\mathbb{N}}$ are uniformly bounded (i.e bounded from above by a constant not depending on $n$). By Lemma \ref{lema} (ii), such a situation is not possible.
\end{rk}

In the estimates of the following proposition there appears an exponent $\delta_0$, which exact value is not known (although according to the Ramanujan conjecture $\delta_0=\frac{1}{2}$). This constant comes from~\cite{AS} and, as it is mentioned there, the best results up to the date are $\frac{25}{64}\leq\delta_0\leq\frac{1}{2}$.

\begin{prop} \label{E}
Let $\alpha$ be a quadratic irrational number, $\D=(u_n)_{n\in\N}$ be a sequence of integers satisfying Hypothesis~\ref{H1}. Fix a sufficiently large index $n\in\N$, more precisely satisfying $u_n>c_\alpha^{-1}$. Denote by $l=l^{(n)}$ the length of the period of $u_n\alpha$ and by $b_1,\ldots,b_l$ the partial quotients in the period of $u_n\alpha$. We have
\begin{equation} \label{proposition_gm_conclusion_ub}
\left(\prod_{i=1}^l b_i\right)\leq e^{\kappa(\gamma+\delta)l},
\end{equation}
where $\gamma:=\frac{\pi^2}{12\log(2)}$ is a so called L\'evy's constant, $\kappa=\frac{24}{\delta_0}$ and $\delta\.=u_n^{-\delta_0/24}$ (the symbol $\dot{=}$ means equal up to a constant independent of $n$). Moreover, 
\begin{equation} \label{proposition_gm_conclusion}
e^{(\gamma-\delta)l}\leq\left(\prod_{i=1}^l (b_i+1)\right)\leq 2^l e^{\kappa(\gamma+\delta)l}.
\end{equation}
\end{prop}

\textbf{Proof.}
Let $x_1=[0,\overline{b_1,\ldots,b_l}]$ and for $2\leq i\leq l$, let $x_i:=[0,\overline{b_i,\ldots,b_l,\ldots, b_{i-1}}]$ be the $(i-1)$-shift of $x_1$.

Define $f:(0,1)\rightarrow\mathbb{R}$ by
$$
f(x)=\begin{cases}
\log(u_n)\delta_0/12& \text{ if } x\in(0,u_n^{-\delta_0/12})\\
-\log(x)& \text{ if } x\in [u_n^{-\delta_0/12},1).
\end{cases}
$$
Note that the function $f(x)$ is $u_n^{\delta_0/12}$-Lipschitz and satisfies $\|f\|_{\infty}=\log(u_n)\delta_0/12$.
Theorem~2.8 from~\cite{AS} applied with $\epsilon=\frac{\delta_0}{24}$ and $q=u_n$ (and taking into account Hypothesis~\ref{H1}) guarantees that
\begin{equation} \label{ie_two}
\left|\int_{0}^{1}f d\nu_{Gauss}\, -\, \dfrac{1}{l}\sum_{i=1}^l f(x_i)\right|\, \ll_{\alpha,S,\delta_0}\, u_n^{-\delta_0/24},
\end{equation}
where $\nu_{Gauss}$ denotes the Gauss-Kuzmin measure  on the unit interval $[0,1]$, that is a measure equivalent to the Lebesgue measure $\nu_{Lebesgue}$ with the relative density $d\nu_{Gauss}=\frac{1}{\log(2)(1+x)} d\nu_{Lebesgue}$.

The integral in the left hand side of~\eqref{ie_two} can be estimated using the following equality:
\begin{equation}
\left|\int_{0}^{1}f(x) d\nu_{Gauss}\right|=\left|\int_{0}^{1}-\log(x) d\nu_{Gauss}+\int_0^{u_n^{-\frac{\delta_0}{12}}}\Big(\frac{\log(u_n)\delta_0}{12}+\log(x)\Big)d\nu_{Gauss}\right|.
\end{equation}
Indeed, on the one hand
\begin{multline}
\left|\int_{0}^{1}\log(x) d\nu_{Gauss}\right|=\left|\int_{0}^{1}\frac{\log(x)}{\log(2)(1+x)} d\nu_{Lebesgue}\right|
\\=\left|\frac{\mathrm{Li}_2(-1)}{\log(2)}\right|=\frac{\pi^2}{12\log(2)}\approx 1.18567,
\end{multline}
where $\nu_{Lebesgue}$ denotes the Lebesgue measure on $\mathbb{R}$ and $\mathrm{Li}_2(z)=\sum_{k=1}^{\infty}\frac{z^k}{k^2}$ is a so called dilogarithm function. On the other hand,
\begin{equation}
\left|\int_0^{u_n^{-\frac{\delta_0}{12}}}\Big(\frac{\log(u_n)\delta_0}{12}+\log(x)\Big)d\nu_{Gauss}\right|\leq \frac{u_n^{-\frac{\delta_0}{12}}}{\log(2)}.
\end{equation}

Thus~\eqref{ie_two} implies
\begin{equation} \label{proposition_gm_db_delta_counting_measure}
\gamma-\delta\,\leq\,\dfrac{1}{l}\sum_{i=1}^l f(x_i)\,\leq\, \gamma+\delta,
\end{equation}
where $\gamma=\frac{\pi^2}{12\log(2)}$ and $\delta\.= u_n^{-\delta_0/24}$.

Since
$$
f(x_i)=\min\left(-\log(x_i),\log(u_n)\delta_0/12\right)=\min\left(\log\Big(\dfrac{1}{x_i}\Big),\log(u_n)\delta_0/12\right)
$$
and $\frac{1}{x_i}<b_i+1$, we have that $f(x_i)\leq \log (b_i+1)$. Thus we infer the lower bound~\eqref{proposition_gm_conclusion} from the lower bound in~\eqref{proposition_gm_db_delta_counting_measure}.

To obtain the upper bound~\eqref{proposition_gm_conclusion}, we split the argument in two cases, according to the two cases in the definition of $f(x)$. Since $b_i<\frac{1}{x_i}$, if $x_i\in [u_n^{-\delta_0/12},1)$, we have
\begin{equation}  \label{proposition_gm_a_i_ub_simple_case}
\log(b_i)<f(x_i)
\end{equation}
(because in this case $f(x_i)=-\log(x_i)$).

In the complimentary case, when $x_i\in (0,u_n^{-\delta_0/12})$, Lemma \ref{lema} (i) gives
$$
b_i\leq \dfrac{u_n}{c_\alpha},
$$
hence
\begin{equation}\label{caso2prop}
\log(b_i)<\frac{12}{\delta_0}f(x_i)-\log(c_\alpha)<\frac{24}{\delta_0}f(x_i)
\end{equation}
(because $u_n>c_\alpha^{-1}$).

We deduce~\eqref{proposition_gm_conclusion_ub} with $\kappa=\frac{24}{\delta_0}$ 
from \eqref{proposition_gm_a_i_ub_simple_case}, \eqref{caso2prop} and the upper bound in~\eqref{proposition_gm_db_delta_counting_measure}.

\begin{flushright}
$\square$
\end{flushright}

\textbf{Proof of Theorem \ref{BZ}.} Let $\alpha$ be a quadratic irrational number and $\D=(u_n)_{n\in\N}$ a sequence satisfying Hypothesis~\ref{H1}. Recall that we want to prove that there exists $c>0$ such that for all $q\in \s'$, where $\s'$ is defined by~\eqref{def_s_prime},
\begin{equation}\label{ec1}
q \cdot \log q \cdot \|q\alpha\| \cdot |q|_{\D}\leq c.
\end{equation}
The inequality \eqref{ec1} will follow from the stronger asymptotic equality
\begin{equation}\label{ec2}
r^{(n)}_{l^{(n)}-1} \cdot \left(\log(r^{(n)}_{l^{(n)}-1}) + \log(u_n)\right) \cdot \|r^{(n)}_{l^{(n)}-1} u_n\alpha\| \,\asymp 1.
\end{equation}
Just by the definition of $\s'$, it is clear that \eqref{ec2} implies \eqref{ec1}. In the rest of the proof we focus on proving~\eqref{ec2}.

We can assume, without loss of generality, that the continued fraction of $\alpha$ is purely periodic, i.e. $\alpha >1$ and its Galois conjugate $\overline{\alpha}$ belongs to the interval $(-1,0)$. Thus the continued fraction expansion of $u_n\alpha$ is almost purely periodic:
$$
u_n\alpha\=[a^{(n)}_0,\overline{a^{(n)}_1,\ldots,a^{(n)}_{l^{(n)}}}].
$$
From the theory of continued fractions, we have
\begin{equation}\label{ec3}
r^{(n)}_{l^{(n)}-1} \|r^{(n)}_{l^{(n)}-1} u_n\alpha\|\asymp\dfrac{1}{a^{(n)}_{l^{(n)}}},
\end{equation}
so, in order to prove \eqref{ec2}, it is enough to show
\begin{equation}\label{ec4}
\log(r^{(n)}_{l^{(n)}-1}) + \log(u_n) \asymp a^{(n)}_{l^{(n)}}.
\end{equation}
To this end, note that from the recurrence satisfied by the denominators of the convergents of $u_n\alpha$,
\begin{equation}\label{rec}
r^{(n)}_0=1,\qquad r^{(n)}_1=a^{(n)}_1,\qquad r^{(n)}_k\=a^{(n)}_k r^{(n)}_{k-1} + r^{(n)}_{k-2}\qquad (k\geq 2),
\end{equation}
we deduce that 
\begin{equation} \label{recurrence_r}
r^{(n)}_k\leq(a^{(n)}_k+1) r^{(n)}_{k-1}
\end{equation}
for all $k\in\N$.

Applying~\eqref{recurrence_r} recursively for $k=1,\dots, l^{(n)}-1$ we find
$$ 
r^{(n)}_{l^{(n)}-1}\leq\prod^{l^{(n)}-1}_{k=0} (a^{(n)}_k+1)
$$
and so, by Proposition \ref{E}, we have 
\begin{equation}\label{ec20}
\log(r^{(n)}_{l^{(n)}-1}) \ll_\alpha l^{(n)}.
\end{equation}
On the other hand, by comparing the recurrence \eqref{rec} with the Fibonacci series (giving the denominators of the convergents of the golden ratio $[\overline{1}]$), we have
\begin{equation}\label{ec7}
\log(r^{(n)}_{l^{(n)}-1})\gg l^{(n)}.
\end{equation}
The inequalities~\eqref{ec20} and~\eqref{ec7} prove
\begin{equation} \label{equiv_r_l}
\log (r_{l^{(n)}-1}^{(n)})\asymp l^{(n)}.
\end{equation}
In \cite{AS}, Theorem 2.12 (1), Aka and Shapira show that 
$l^{(n)}\asymp u_n$. Then, applying this result to~\eqref{equiv_r_l} we obtain
\begin{equation} \label{ec21}
\log(r^{(n)}_{l^{(n)}-1})\asymp u_n.
\end{equation}
By a trivial asymptotic comparison we infer from~\eqref{ec21}
$$
\log(r^{(n)}_{l^{(n)}-1})+\log(u_n)\asymp u_n.
$$
Finally, by Lemma \ref{lema} we have $u_n\asymp a^{(n)}_{l^{(n)}}$, thus
$$
\log(r^{(n)}_{l^{(n)}-1})+\log(u_n)\asymp a^{(n)}_{l^{(n)}}.
$$
This establishes~\eqref{ec4} and concludes the proof of the theorem.
\begin{flushright}
$\square$
\end{flushright}

\end{document}